\input amstex
%\input epsf
%\input preamble
%\newsymbol\boxtimes 1202

%letters

\def\b1{\text{\bf 1}}

\def\BC{{\Bbb C}}

\def\BZ{{\Bbb Z}}
\def\CA{{\Cal A}}

\def\CD{{\Cal D}}

\def\CL{{\Cal L}}
\def\CM{{\Cal M}}

\def\CO{{\Cal O}}

\def\CT{{\Cal T}}

\def\CV{{\Cal V}}

\def\dpar{\partial}
\def\dplus{\buildrel\cdot\over{+}}

\def\fa{{\frak a}}

\def\fb{{\frak b}}

\def\fD{{\frak D}}

\def\fg{{\frak g}}

\def\fh{{\frak h}}

\def\fn{{\frak n}}

\def\fp{{\frak p}}

\def\hfg{\widehat{\frak g}}

\def\htau{\widehat{\tau}}

\def\tC{\tilde C}
\def\tCA{\widetilde{\Cal{A}}}

\def\tCV{\widetilde{\Cal{V}}}

\def\tf{\tilde f}

%symbols

\def\btu{\bigtriangleup}
\def\hra{\hookrightarrow}
\def\iso{\buildrel\sim\over\longrightarrow} 

\def\lra{\longrightarrow}

% space between paragraphs 
\parskip=6pt

\documentstyle{amsppt}
\document
%\magnification=1200
\NoBlackBoxes
%\nologo

%%%%%%%%%

\centerline{\bf On chiral differential operators over homogeneous spaces}

\bigskip
\centerline{Vassily Gorbounov, Fyodor Malikov, Vadim Schechtman}
\bigskip

%\centerline{\it (Preliminary draft)} 

%\bigskip\bigskip

%\input groupintro
\centerline{\bf Introduction}

\bigskip\bigskip 

The notion of an {\it algebra of chiral differential operators} 
(cdo for short) over a smooth 
algebraic variety $X$ has been studied in [GII]. (This notion 
has been invented and first studied, in a different language, by 
Beilinson and Drinfeld, cf. [BD1], Chapter 3, \S 8).  
In the present note we consider some examples in more detail. 
We shall work over the ground field $\BC$. 

First, we give a classification of cdo over $X$ in the following cases: 
$X=G$ is  an affine algebraic group; $X=G/N$ or $G/P$ where $N$ is a unipotent 
subgroup and $P$ is a parabolic subgroup and $G$ is simple (the extension 
to the case of a semisimple $G$ being straightforward). 

Before we describe the result, let us explain some terminology and 
notation. For a smooth algebraic variety $X$, an algebra of cdo 
over $X$ is by definition a Zariski sheaf $\CV$ of $\BZ_{\geq 0}$-graded 
vertex algebras on $X$ such that 

(a) if $Alg(\CV)=(\CA,\CT,\Omega,\dpar,\gamma, 
\langle,\rangle,c)$ is the sheaf of {\it vertex algebroids} associated 
with $\CV$ 
(see [GII], $\S 2$), then the corresponding {\it extended Lie algebroid}  
$(\CA,\CT,\Omega,\dpar)$ ({\it op. cit.}, 1.1, 1.4) is identified with 
$(\CO_X,\Theta_X,\Omega^1_X,d_{DR})$ where $\Theta_X$ denotes the tangent 
bundle and $d_{DR}$ the de Rham differential;  

(b) the adjunction morphism $U Alg(\CV)\lra\CV$ is an isomorphism. Here 
$U$ is the functor of {\it vertex envelope} defined in {\it op. cit.}, 
\S 9. 

For each Zariski open $U\subset X$ we can consider the category 
(a groupoid in fact) of 
cdo over $U$, or, what is the same, the groupoid 
of vertex algebroids (defined in [GII], \S 3)  
over $U$ satisfying (a) above. When $U$ varies, we get a sheaf of 
groupoids $\fD iff^{ch}_X$ over $X$ --- the {\it gerbe of chiral 
differential operators}. As usual, $\Gamma(U;\fD iff^{ch}_X)$ will denote 
the sections over $U$; a generic object of this category will 
be sometimes denoted $\CD^{ch}_U$; the set of isomorphism classes 
of cdo over $U$ will be denoted by 
$\pi_0(\Gamma(U;\fD iff^{ch}_X))$. 

Let $G$ be an affine connected algebraic group, $\fg$ the 
corresponding Lie algebra. For each symmetric ad-invariant bilinear form 
$(,)\in (S^2\fg^*)^{\fg}$ we construct a cdo $\CD^{ch}_{G;(,)}$ 
over $G$ such that if $G$ is semisimple, then the correspondence 
$(,)\mapsto\CD^{ch}_{G;(,)}$ gives rise to a bijection 
$$
(S^2\fg^*)^{\fg}\iso\pi_0(\Gamma(G,\fD iff^{ch}_G))
\eqno{(I1)}
$$          
We have a canonical embedding of vertex algebras 
$$
i_{(,)}:\ \CV_{\fg,(,)}\hra \CD^{ch}_{G;(,)}
\eqno{(I2)}
$$
where $\CV_{\fg;(,)}$ denotes the vacuum module of the Kac-Moody algebra 
$\hfg$ at level $(,)$. This embedding is induced by the 
embedding of $\fg$ into $T_G:=\Gamma(G;\Theta_G)$ as left invariant 
vector fields. 

Let $(,)_{\fg;(K)}$ denote the Killing form on $\fg$. Let us define 
the {\it dual} level by 
$$
(,)^o=-(,)_{\fg;(K)}-(,)
\eqno{(I3)}
$$ 
Using the embedding $\fg\hra T_G$ by means of {\it right} invariant 
vector fields one can construct a canonical {\it dual embedding}  
of vertex algebras 
$$
i^o_{(,)}:\ \CV_{\fg;(,)^o}\hra\CD^{ch}_{G;(,)}
\eqno{(I4)}
$$
It is characterised by the requirement that the images of $i_{(,)}$ 
and $i^o_{(,)}$ {\it commute} in an appropriate sense (see Theorem 2.5 
and Corollary 2.6). This beautiful fact was communicated to us 
by B.Feigin, E.Frenkel and D.Gaitsgory. We give a proof using the language 
of [GII]. 

Let us pass to homogeneous spaces. Assume that $G$ is simple. 
Let $N\subset G$ be a unipotent group. The classification of cdo over 
$G/N$ is the same as over $G$; namely, for each level $(,)$ 
one can define a cdo $\CD^{ch}_{G/N;(,)}$ such that the correspondence 
$(,)\mapsto \CD^{ch}_{G/N;(,)}$ induces a bijection 
$$
(S^2\fg^*)^{\fg}\iso\pi_0(\Gamma(G/N;\fD iff^{ch}_{G/N}))
\eqno{(I5)}
$$   
The sheaves $\CD^{ch}_{G/N;(,)}$ are constructed using the {\it BRST}  
(or quantum Hamiltonian) reduction of the corresponding cdo's on $G$. 
More precisely, 
$$
\CD^{ch}_{G/N; (,)}=H_{BRST}^{0}(L\fn;\pi_*\CD^{ch}_{G;(,)})
\eqno{(I6)}
$$
where the rhs denotes the BRST cohomology of the loop algebra 
$L\fn:=\fn [T,T^{-1}],\ \fn:=Lie(N)$. For the precise definition see 
\S 3. 

Let $B\subset G$ be a Borel subgroup. We show that there exists a unique, 
up to a unique isomorphism, cdo $\CD^{ch}_{G/B}$ on the flag space $G/B$. 
Again this cdo may be constructed using the BRST reduction. Namely, 
$$
\CD^{ch}_{G/B}=H_{BRST}^{0}(L\fb,\fh;\pi_*\CD^{ch}_{G;crit})
\eqno{(I7)}
$$
Here $\CD^{ch}_{G;crit}$ is by definition the cdo $\CD^{ch}_{G;(,)_{crit}}$ 
on the {\it critical} level $(,)_{crit}=-(,)_{\fg;(K)}/2$. 
For the definition of the relative BRST cohomology in the 
rhs we again refer the reader to the main body of the note, see \S 4. 
(A more explicit construction of the sheaf 
$\CD^{ch}_{G/B}$ for $G=SL(n)$, using vertex operators, has been suggested 
in [MSV], 5.9-5.10.)  

The embeddings (I4) induce canonical morphisms of vertex algebras 
$$
\CV_{\fg;(,)^o}\lra \CD^{ch}_{G/N;(,)};\ \CV_{\fg;(,)_{crit}}\lra 
\CD^{ch}_{G/B}
\eqno{(I8)}
$$       
Taking the spaces of sections over a big cell, we get another construction 
of Feigin-Frenkel {\it Wakimoto modules}, cf. [FF1] - [FF3]. 

Finally, if $P\subset G$ is  parabolic but not Borel, we show that 
$\Gamma(G/P;\fD iff^{ch}_{G/P})$ is empty. The classification 
of cdo over homogeneous spaces is exactly reflected in the BRST world: 
namely the square of the corresponding BRST charge is zero at all 
levels for $G/N$, only at the critical level for $G/B$ and 
is never zero for $G/P$.  

This introduction would not be complete without mentioning that 
this note relies heavily on the ideas of B.Feigin and E.Frenkel. 
This note started   
from our attempts to find a proof of 2.5 and 2.6.
Our sincere gratitude goes to D.Gaitsgory who had communicated these facts
to us and told us that he had known their proofs. 
We are also grateful to H.Esnault for a crucial remark 4.1.1.  

This work was done while the authors were visiting Institut des Hautes 
\'Etudes Scientifiques in Bures-sur-Yvette and 
Max-Planck-Institut f\"ur Mathematik in Bonn. We are grateful 
to these institutes for the hospitality.

\bigskip\bigskip

\centerline{\bf \S 1. Chiral differential operators over an 
algebraic group} 

\bigskip\bigskip  

{\it Perfect vertex algebroids over constants}

\bigskip

The discussion below is nothing but the specification of [GII], \S\S 1 - 4 
to the case $A=\BC$. 

{\bf 1.1.} Let $\fg$ be a Lie algebra. We shall need 
two complexes connected with $\fg$, both concentrated in nonnegative 
degrees. The first one, $C^{\cdot}(\fg)=
C^{\cdot}(\fg;\BC)$, is the cochain complex of $\fg$ with trivial 
coefficients. Thus, by definition 

$C^i(g)=(\Lambda^i\fg)^*=$ the space of skew symmetric polylinear maps 
$f:\ \fg^i\lra \BC$, $i\geq 0$.  

The differential $d:\ C^{i-1}(\fg)\lra C^i(\fg)$ acts as 
$$
df(\tau_1,\ldots,\tau_i)=\sum_{1\leq p<q\leq i}\ 
(-1)^{p+q+1}f([\tau_p,\tau_q],\tau_1,\ldots,\htau_p,\ldots,\htau_q,
\ldots,\tau_i)
\eqno{(1.1.1)}
$$
The cohomology spaces $H^i(C^{\cdot}(\fg))$ will be denoted $H^i(\fg)$. 
   
The second complex, $\tC^{\cdot}(\fg)$, is the shifted by $1$ and augmented  
cochain complex of $\fg$ with coefficients in the coadjoint representation 
$\fg^*$. By definition, $\tC^0(\fg)=\BC$ and 

$\tC^i(\fg)=Hom_k(\Lambda^{i-1}\fg,\fg^*)=$ the space of skew symmetric 
polylinear maps 
\newline $h:\ \fg^{i-1}\lra\fg^*$ for $i\geq 1$. 

The differential $d:\ \tC^0(\fg)\lra\tC^1(\fg)$ is zero, and 
$d:\ \tC^i(\fg)\lra\tC^{i+1}(\fg)$ acts as 
$$
dh(\tau_1,\ldots,\tau_i)=\sum_{p=1}^i\ (-1)^p\tau_p(h(\tau_1,\ldots,
\htau_p,\ldots,\tau_i))+
$$
$$
+\sum_{1\leq p<q\leq i}\ 
(-1)^{p+q}h([\tau_p,\tau_q],\tau_1,\ldots,\htau_p,\ldots,\htau_q,
\ldots,\tau_i)
\eqno{(1.1.2)}
$$
for $i\geq 1$. Define embeddings $C^i(\fg)\hra\tC^i(\fg)$ by assigning 
to $f\in C^i(\fg)$ an element $\tf\in\tC^i(\fg)$ given by 
$$
\langle\tau_1,\tf(\tau_2,\ldots,\tau_i)\rangle=
f(\tau_1,\ldots,\tau_i)
\eqno{(1.1.3)}
$$
We shall identify $C^i(\fg)$ with its image in $\tC^i(\fg)$. 

One checks that the embeddings (1.1.3) are compatible with the differentials, 
so that one has an embedding of complexes $C^{\cdot}(\fg)\hra
\tC^{\cdot}(\fg)$. 

{\bf 1.2.} Let us consider the groupod $\CA lg_{\fg}$ of vertex algebroids 
of the form 
\newline $\CA=(\BC,\fg,\fg^*,\dpar,\gamma,\langle,\rangle,c)$ where 
$\CT=(\BC,T,\Omega,\dpar)=(\BC,\fg,\fg^*,0)$ is a {\it perfect} 
extended Lie algebroid 
over $\BC$, in the sense of [GII], 1.2, with $T=\fg$. Note that the 
last object is uniquely 
defined by the Lie algebra $\fg=T$; we must have $\Omega=\fg^*$, the 
"Lie derivative" action of $T$ on $\Omega$ must be the coadjoint one, 
and a $\BC$-linear derivation $\dpar:\ \BC\lra\Omega$ must be zero. 

Turning to the axioms of a vertex algebroid, {\it op. cit.}, 1.4, we see 
that for $\CA$ as above, $\langle,\rangle:\ \fg\times\fg\lra\BC$ is a symmetric 
bilinear map (which may be regarded as an element of $\tC^2(\fg)$), 
$c\in\tC^3(\fg)$,  
(A1) implies that $\gamma=0$, (A2) and (A3) 
hold true automatically, (A4) takes the form 
$$
\langle [\tau_1,\tau_2],\tau_3\rangle+
\langle \tau_2,[\tau_1,\tau_3]\rangle=\langle\tau_2,c(\tau_1,\tau_3)\rangle+
\langle\tau_3,c(\tau_1,\tau_2)\rangle
\eqno{(1.2.1)}
$$
and (A5) takes the form 
$$
dc=0
\eqno{(1.2.2)}
$$
where $d$ is the differential in $\tC^{\cdot}(\fg)$ given by (1.1.2). 

So, an object of $\CA lg_{\fg}$ has a form 
$$
\CA_{\fg;\langle\rangle,c}=(\BC,\fg,\fg^*,0,0,\langle,\rangle,c)
\eqno{(1.2.3)}
$$ 
where $\langle,\rangle\in\tC^2(\fg)^{\BZ/2\BZ},\ 
c\in\tC^3(\fg)$ satisfy (1.2.1) and (1.2.2). 

The vertex envelope 
$$
\CV_{\fg;\langle,\rangle,c}=U\CA_{\fg;\langle,\rangle,c}
\eqno{(1.2.4)}
$$
(see [GII], \S 9) is generated by the fields $\tau(z)\ (\tau\in\fg)$ and 
$\omega(z)\ (\omega\in\fg^*)$ of conformal weight $1$, subject to OPE 
$$
\tau(z)\tau'(w)\sim\frac{\langle\tau,\tau'\rangle}{(z-w)^2}+
\frac{[\tau,\tau'](w)-c(\tau,\tau')(w)}{z-w}
\eqno{(1.2.5)}
$$
$$
\tau(z)\omega(w)\sim\frac{\tau(\omega)(w)}{z-w};\ 
\omega(z)\omega'(w)\sim 0
\eqno{(1.2.6)}
$$
cf. [GII] (9.9.1)--- (9.9.3).  

A morphism 
$$
f:\ \CA_{\fg,\langle,\rangle,c}\lra 
\CA_{\fg,\langle,\rangle',c'}
$$
is by definition an element $h\in \tC^2(\fg)$ such that 
$$
\langle\tau_1,h(\tau_2)\rangle+\langle\tau_2,h(\tau_1)\rangle=
\langle\tau_1,\tau_2\rangle- \langle\tau_1,\tau_2\rangle'
\eqno{(1.2.7)}
$$
and 
$$
dh=c-c'
\eqno{(1.2.8)}
$$
see [GII], Theorem 3.5. The composition of morphisms is induced by the 
addition in $\tC^2(\fg)$.  

{\bf 1.3.} As a corollary, we have a canonical bijection 
$$
\pi_0(\CA lg_{\fg})=H^3(\fg)
\eqno{(1.3.1)}
$$
More precisely, for a $3$-cocycle $c\in C^{3,cl}(\fg)$ we have  
a vertex algebroid 
$$
\CA_{\fg;c}:=\CA_{\fg;0,c}
\eqno{(1.3.2)}
$$
and the correspondence $c\mapsto \CA_{\fg;c}$ induces the bijection (1.3.1). 

The enveloping algebra $\CV_{\fg;c}:=U\CA_{\fg;c}$ is generated by the same 
fields as in 1.2, subject to OPE
$$
\tau(z)\tau'(w)\sim 
\frac{[\tau,\tau'](w)-c(\tau,\tau')(w)}{z-w}
\eqno{(1.3.3)}
$$
and (1.2.6).   

Let us define another interesting class of objects of $\CA lg_{\fg}$. 
Namely, each symmetric $ad$-invariant bilinear form $(,)\in (S^2\fg^*)^{\fg}$ 
gives rise to an object 
$$
\tCA_{\fg;(,)}:=\CA_{\fg;(,),0}
\eqno{(1.3.4)}
$$
The enveloping algebra $\CV_{\fg;c}:=U\CA_{\fg;c}$ is generated by the same 
fields as in 1.2, subject to OPE
$$
\tau(z)\tau'(w)\sim\frac{(\tau,\tau')}{(z-w)^2}+
\frac{[\tau,\tau'](w)}{z-w}
\eqno{(1.3.5)}
$$
and (1.2.6). 

It is easy to see that we have an isomorphism 
$$
f_{(,)}:\ \tCA_{\fg;(,)}\iso \CA_{\fg;c_{(,)}}
\eqno{(1.3.6)}
$$
given by a map $h_{(,)}:\fg\lra\fg^*$ where
$$
\langle\tau_1,h_{(,)}(\tau_2)\rangle=\frac{1}{2}(\tau_1,\tau_2)
\eqno{(1.3.7)}
$$
and the cocycle $c_{(,)}$ is defined by 
$$
c_{(,)}(\tau_1,\tau_2,\tau_3)=([\tau_1,\tau_2],\tau_3)
\eqno{(1.3.8)}
$$
Cf. [GII], Theorem 4.5. 

{\bf 1.4.} Given $(,)\in (S^2\fg^*)^{\fg}$, consider a vertex algebroid 
$$
\CA_{\fg;(,)}:=(\BC,\fg,\fg^*,0,0,(,),0)
\eqno{(1.4.1)}
$$
Its vertex envelope $\CV_{\fg;(,)}:=U\CA_{\fg;(,)}$ is generated by fields 
$\tau(z)\ (\tau\in\fg)$ of conformal weight $1$, subject to OPE (1.3.5). 

The correspondence $\tau\cdot T^n\mapsto \tau_{(n)}$ defines on 
$\CV_{\fg;(,)}$ a structure of the vacuum module over the Kac-Moody  
algebra $\hfg=\fg[T, T^{-1}]\oplus\BC\cdot\b1$ at level $(,)$. 

We have an obvious embedding of vertex algebroids 
$\CA_{\fg;(,)}\hra\CA_{\fg;0}$ which induces an embedding of vertex 
algebras 
$$
\CV_{\fg;(,)}\hra\tCV_{\fg;(,)}
\eqno{(1.4.2)}
$$   

{\bf 1.5.} If the Lie algebra $\fg$ is {\it semisimple} then the 
correspondence $(,)\mapsto c_{(,)}$ induces a bijection 
$$
(S^2\fg^*)^{\fg}\iso H^3(\fg)
\eqno{(1.5.1)}
$$  
Therefore, in this case the algebroids $\tCA_{\fg;(,)}$ form a complete 
set of representatives of isomorphism classes in $\CA lg_{\fg}$. In other 
words, 

{\bf 1.5.1.} {\it if $\fg$ is semisimple then the correspondence 
$(,)\mapsto\tCA_{\fg;(,)}$ induces a bijection 
$$
(S^2\fg^*)^{\fg}\iso \pi_0(\CA lg_{\fg})
\eqno{(1.5.1.1)}
$$}

\bigskip

{\it Passing to a group}

\bigskip

{\bf 1.6.} Let $G=Spec(A)$ be an affine algebraic group, 
$\fg$  the corresponding Lie algebra . The tangent bundle $\Theta_G$ is trivial, 
so the obstruction 
$c(\fD iff^{ch}_G)$ to the existence of an algebra of chiral do over $G$, 
$\CD^{ch}_G\in\Gamma(G;\fD iff^{ch}_G)$ (cf. [GII], Corollary 7.11) vanishes. 

By {\it op. cit.}, \S 4, the set of isomorphism classes of chiral do 
over $G$, $\pi_0(\Gamma(G;\fD iff^{ch}_G))$ is a nonempty torseur under 
the "Chern-Simons group" $H^3_{DR}(G)=H^3(G;\BC)$. 

In fact the groupoid $\Gamma(G;\fD iff^{ch}_G)$ has a distinguished object 
$\CD^{ch}_{G;0}$, so that we have a canonical bijection 
$$
\pi_0(\Gamma(G;\fD iff^{ch}_G))\iso H^3_{DR}(G)
\eqno{(1.6.1)}
$$            
This is a consequence of the following general construction. 

{\bf 1.7.} Let $\CA_{\fg;\langle,\rangle,c}$ be an arbitrary object 
of $\CA lg_{\fg}$. Let us apply to it the {\it pushout} construction of 
[GII], 1.10 with respect to the structre morphism $\BC\lra A$. 
Here the morphism $\fg\lra T:=Der_{\BC}(A)$ is defined as the embedding of 
{\it left invariant} vector fields, and the map 
$\gamma:\ A\times\fg\lra \Omega:=\Omega^1(A)$ is set to be zero. 
This way we get a vertex $A$-algebroid $\CA_{G;\langle,\rangle,c}$. 
Its enveloping algebra
$$
\CD^{ch}_{G;\langle,\rangle,c}=U\CA_{G;\langle,\rangle,c}
\eqno{(1.7.1)}
$$
obviously belongs to $\Gamma(G;\fD iff^{ch}_G)$. 

We have a canonical embedding 
$$
\CV_{\fg,\langle,\rangle,c}\hra\CD^{ch}_{G;\langle,\rangle,c}
\eqno{(1.7.2)}
$$ 

We shall use the notations $\CA_{G;(,)}:=\CA_{G;(,),0},\ 
\CA_{G;c}:=\CA_{G;0,c},\ \CA_{G;0}:=\CA_{G;0,0}$ and 
$\CD^{ch}_{G;(,)}$, etc. for the corresponding enveloping algebras. 

If $\CA_{G;0}=(A,T,\Omega,d_{DR},\gamma_0,\langle,\rangle_0,c_0)$ and 
$\omega\in \Omega^{3,cl}(A)$ is a closed $3$-form then we can form 
a vertex algebroid 
$$
\CA_{G;\omega}:=\CA_{G;0}\dplus\omega=
(A,T,\Omega,d_{DR},\gamma_0,\langle,\rangle_0,c_0+\omega)
\eqno{(1.7.3)}
$$   
The correspondence $\omega\mapsto\CA_{G;\omega}$ induces the bijection 
(1.6.1). 

If $c\in C^{3,cl}(\fg)$ is a $3$-cocycle with trivial coefficients then 
by definition  
$$
\CA_{G;c}=\CA_{G;\omega_c}
\eqno{(1.7.4)}
$$
where $\omega_c\in\Omega^{3,cl}(A)$ is the left invariant $3$-form on $G$ 
corresponding to $c$. 

{\bf 1.8. Corollary.} {\it Assume that $G$ is reductive. Then 
the correspondence $c\mapsto\CA_{G;c}$ induces a bijection 
$$
H^3(\fg)\iso \pi_0(\Gamma(G;\fD iff^{ch}_G))
\eqno{(1.8.1)}
$$}

Indeed, one knows that for a reductive group the correspondence 
$c\mapsto\omega_c$ gives rise to an isomorphism $H^3(\fg)\iso H^3_{DR}(G)$. 

{\bf 1.9. Corollary.} {\it Assume that $G$ is semisimple. Then 
the correspondence $(,)\mapsto\CA_{G;(,)}$ induces a bijection 
$$
(S^2\fg^*)^{\fg}\iso \pi_0(\Gamma(G;\fD iff^{ch}_G))
\eqno{(1.9.1)}
$$}

This follows from 1.8 and the remarks 1.5. 

{\bf 1.10.} Note that for an arbitrary $G$ and 
$(,)\in (S^2\fg^*)^{\fg}$ one has a canonical embedding 
$$
i_{(,)}:\ \CV_{\fg;(,)}\hra\CD^{ch}_{G;(,)}
\eqno{(1.10.1)}
$$
It is the composition of (1.4.2) and (1.7.2).       

\bigskip\bigskip

\newpage

\centerline{\bf \S 2. Dual embedding}

\bigskip\bigskip

{\bf 2.1.} Let $G=Spec(A)$ be a smooth affine connected algebraic group with 
the Lie algebra $\fg$. Pick a symmetric ad-invariant bilinear form 
({\it "level"}) $(,)\in (S^2\fg^*)^{\fg}$. 

Let $(,)_{(K)}$ denote the Killing form on $\fg$, 
$$
(x,y)_{(K)}=tr_{\fg}(ad_x\cdot ad_y)
\eqno{(2.1.1)}
$$
Let us pick a base $\{\tau_i\}$ of $\fg$. In terms of structure 
constants
$$
[\tau_i,\tau_j]=c^{ij}_p\tau_p
\eqno{(2.1.2)}
$$
the form (2.1.1) is given by 
$$
(\tau_i,\tau_j)_{(K)}=c^{ip}_qc^{jq}_p
\eqno{(2.1.3)}
$$
Let us define the {\it dual level} $(,)^o\in (S^2\fg^*)^{\fg}$ by 
$$
(,)^o=-(,)_{(K)}-(,)
\eqno{(2.1.4)}
$$
Define the {\it critical level} $(,)_{crit}$ by $(,)_{crit}=(,)^o_{crit}$, 
i.e. 
$$
(,)_{crit}=-\frac{1}{2}(,)_{(K)}
\eqno{(2.1.5)}
$$
If we want to stress the dependence on $\fg$, we shall write 
$(,)_{\fg;(K)},\ (,)_{\fg; crit}$.  

{\bf 2.2.} We have two commuting left actions of $G$ on itself: 
the left multiplication, $(g,x)\mapsto gx$ and the right one, 
$(g,x)\mapsto xg^{-1}$. 

Let $T=Der_k(A)$ denote the Lie $A$-algebroid of vector fields 
over $G$. The above two actions induce two embeddings of Lie algebras 
$$
i_L:\ \fg\hra T,\ i_R:\ \fg\hra T
\eqno{(2.2.1)}
$$
such that 
$$
[i_L(x),i_R(y)]=0,\ \text{for all\ }x,y\in\fg
\eqno{(2.2.2)}
$$
Below we shall identify $\fg$ with its image under $i_L$, i.e. 
write simply $x$ instead of $i_L(x)$. We shall also use the notation 
$x^R:=i_R(x)\ (x\in\fg)$. 

Embedding $i_L$ induces an isomorphism of left $A$-modules 
$$
A\otimes_k\fg\iso T
\eqno{(2.2.3)}
$$
Thus, $\{\tau_i\}$ form an $A$-base of $T$. In particular 
$$
\tau_i^R=a^{ij}\tau_j
\eqno{(2.2.4)}
$$
for some invertible matrix $(a^{ij})$ over $A$. 

The commutation relations 
$$
[\tau_i,\tau_j^R]=0
\eqno{(2.2.5)}
$$ 
are equivalent to the identities 
$$
\tau_i(a^{js})+c^{ip}_sa^{jp}=0
\eqno{(2.2.6)}
$$
true for all $i, j, s$. 

Let us write down the relations 
$$
[\tau_i^R,\tau_j^R]=[\tau_i,\tau_j]^R
\eqno{(2.2.7)}
$$
in coordinates. We have 
$$
[\tau_i^R,\tau_j^R]=[\tau_i^R,a^{js}\tau_s]=\tau_i^R(a^{js})\tau_s=
a^{ip}\tau_p(a^{js})\tau_s
\eqno{(2.2.8)}
$$
due to (2.2.5). Plugging this into (2.2.7) we get 
$$
a^{ip}\tau_p(a^{js})=c^{ij}_qa^{qs}
\eqno{(2.2.9)}
$$
for all $i, j, s$.    

{\bf 2.3.} We set $\Omega:=\Omega^1_{A/k}=Hom_A(T,A)$, and denote by 
$\langle,\rangle:\ T\times\Omega\lra A$ the canonical $A$-bilinear pairing. 
Let $\{\omega_i\}$ be the $A$-base of $\Omega$ dual to $\{\tau_i\}$. 
The Lie algebra $T$ acts on $\Omega$ by the Lie derivative.  

We have $\tau_i(\omega_j)=\alpha^{ijs}\omega_s$ where 
$$
\alpha^{ijs}=\langle\tau_s,\tau_i(\omega_j)\rangle=
\tau_i(\langle\tau_s,\omega_j\rangle)-
\langle [\tau_i,\tau_s],\omega_j\rangle=-c^{is}_j=c^{si}_j
$$
Thus, 
$$
\tau_i(\omega_j)=c^{si}_j\omega_s
\eqno{(2.3.1)}
$$
Similarly, 
$$
\tau_i^R(\omega_j)=0
\eqno{(2.3.2)}
$$

{\bf 2.4.} Recall that we have an embedding of vertex algebras 
$$
i_{(,)}:\ \CV_{\fg;(,)}\hra\CD^{ch}_{G;(,)}
\eqno{(2.4.1)}
$$
see 1.10. More precisely, it is induced by an embedding of conformal 
weight $1$ components
$$
j_L:\ \fg=\CV_{\fg,(,)1}\lra \CD^{ch}_{G;(,)1}=T\oplus\Omega
\eqno{(2.4.2)}
$$
defined by a composition 
$$
\fg\hra T\hra T\oplus\Omega
$$
where the first arrow is $i_L$ and the second one sends $x$ to $(x,0)$. 

The fact that $j_L$ induces a map of vertex algebras (2.4.1) simply means 
that we have the identities in $\CD^{ch}_{G;(,)}$     
$$
j_L(\tau)_{(1)}j_L(\tau')=(\tau,\tau');\ 
j_L(\tau)_{(0)}j_L(\tau')=j_L([\tau,\tau'])
\eqno{(2.4.3)}
$$
for all $\tau, \tau'\in\fg$. 

{\bf 2.5. Theorem.} (B.Feigin - E.Frenkel, D.Gaitsgory) 
(i) {\it There exists a unique embedding 
$$
j_R:\ \fg\hra \CD^{ch}_{G;(,)1}
\eqno{(2.5.1)}
$$
such that}

(a) {\it the composition of} (2.5.1) {\it with the canonical projection 
$\CD^{ch}_{G;(,)1}\lra T$ is equal to $i_R$}; 

(b) {\it for all $\tau, \tau'\in\fg$ and $n\geq 0$
$$
j_L(\tau)_{(n)}j_R(\tau')=0
\eqno{(2.5.2)}
$$} 

(ii) {\it We have 
$$
j_R(\tau)_{(1)}j_R(\tau')=(\tau,\tau')^o
\eqno{(2.5.3)}
$$
and
$$
j_R(\tau)_{(0)}j_R(\tau')=j_R([\tau,\tau'])
\eqno{(2.5.4)}
$$
for each $\tau, \tau'\in\fg$.} 

{\bf 2.6. Corollary.} (B.Feigin - E.Frenkel, D.Gaitsgory) {\it The map} (2.5.1) 
{\it 
induces an embedding of chiral algebras
$$
j_R:\ \CV_{\fg;(,)^o}\hra \CD^{ch}_{G;(,)}
\eqno{(2.6.1)}
$$}

{\it The images of $j_L$ and $j_R$} commute {\it in the following sense: 
$$
j_L(x)_{(n)}j_R(y)=0
\eqno{(2.6.2)}
$$
for each $x\in\CV_{\fg;(,)},\ y\in\CV_{\fg;(,)^o}$ and $n\geq 0$.}

{\bf 2.7.} {\it Proof of} (2.5). (i) As usual, we denote $j_L(\tau)$ 
simply by $\tau$. 
We are looking for $j_R(\tau)$ in the form 
$$
j_R(\tau_i)=\tau_i^R+b^{iq}\omega_q
\eqno{(2.7.1)}
$$
for some $b^{iq}\in A$. We have, by [GII], 1.4 (A2) 
$$
\tau^R_{i(1)}\tau_j=\langle\tau^R_i,\tau_j\rangle=
\langle a^{ip}\tau_p,\tau_j\rangle=
a^{ip}(\tau_p,\tau_j)-\tau_p\tau_j(a^{ip})
$$
Using (2.2.6) and (2.1.4), 
$$
-\tau_p\tau_j(a^{ip})=\tau_p(c^{js}_pa^{is})=
-c^{pu}_pc^{js}_pa^{iu}=c^{up}_sc^{js}_pa^{iu}=
(\tau_u,\tau_j)_{(K)}a^{iu},
$$
so
$$
\tau^R_{i(1)}\tau_j=-(\tau_p,\tau_j)^oa^{ip}
\eqno{(2.7.2)}
$$
On the other hand, 
$$
(b^{iq}\omega_q)_{(1)}\tau_j=\langle b^{iq}\omega_q,\tau_j\rangle=b^{ij}
$$
Therefore the condition 
$$
j_R(\tau_i)_{(1)}\tau_j=0
\eqno{(2.7.3)}
$$
defines the matrix $(b^{iq})$ uniquely: we must have 
$$
j_R(\tau_i)=\tau^R_i+(\tau_p,\tau_q)^oa^{ip}\omega_q
\eqno{(2.7.4)}
$$

{\bf 2.8.} Let us prove that 
$$
\tau_{i(0)}j_R(\tau_j)=0
\eqno{(2.8.1)}
$$
We have 
$$
\tau_{i(0)}j_R(\tau_j)=\tau_{i(0)}\bigl\{\tau_j^R+
(\tau_s,\tau_u)^oa^{js}\omega_u\bigl\}
$$
On the one hand, 
$$
\tau_{i(0)}\tau_j^R=\tau_{i(0)}(a^{jq}\tau_q)=
\tau_{i(0)}(a^{jq}_{(-1)}\tau_q)=
\tau_i(a^{jq})_{(-1)}\tau_q+a^{jq}_{(-1)}[\tau_i,\tau_q]=0
$$ 
using (2.2.6). On the other hand, 
$$
\tau_{i(0)}\bigl\{(\tau_s,\tau_u)^oa^{js}\omega_u\bigr\}=
(\tau_s,\tau_u)^o\tau_i(a^{js}\omega_u)=
(\tau_s,\tau_u)^o\bigl\{\tau_i(a^{js})\omega_u+
a^{js}\tau_i(\omega_u)\bigr\}=
$$
using (2.2.6) and (2.3.1)
$$
=-(\tau_s,\tau_p)^oc^{iq}_sa^{jq}\omega_p-
(\tau_s,\tau_u)^oc^{ip}_ua^{js}\omega_p=
-([\tau_i,\tau_q],\tau_p)^oa^{jq}\omega_p
-(\tau_s,[\tau_i,\tau_p])^oa^{js}\omega_p=0
$$
due to the invariance of the form $(,)^o$. This proves (2.8.1). 

Evidently $\tau_{i(n)}j_R(\tau_j)=0$ for $n\geq 2$. This proves part (i) 
of the theorem. 

{\bf 2.9.} Let us compute $j_R(\tau_i)_{(1)}j_R(\tau_j)$. We have 
$$
j_R(\tau_i)_{(1)}j_R(\tau_j)=\bigl\{\tau_i^R+(\tau_p,\tau_q)^o
a^{ip}\omega_q\bigr\}_{(1)}
\bigl\{\tau_j^R+(\tau_s,\tau_u)^oa^{js}\omega_u\bigr\}
$$   
This is a sum of four terms. 
$$
I:=\tau^R_{i(1)}\tau_j^R=\langle a^{ip}\tau_p,a^{js}\tau_s\rangle=
$$
using [GII] (1.8.3)$_{\langle,\rangle}$
$$
=a^{ip}a^{js}(\tau_p,\tau_s)-
a^{ip}\tau_s\tau_p(a^{js})-a^{js}\tau_p\tau_s(a^{ip})-
\tau_p(a^{js})\tau_s(a^{ip})
$$
Using (2.2.6) and (2.1.4) we see that 
$$
-a^{ip}\tau_s\tau_p(a^{js})=-a^{js}\tau_p\tau_s(a^{ip})=
\tau_p(a^{js})\tau_s(a^{ip})=(\tau_p,\tau_s)_{(K)}a^{ip}a^{js},
$$
whence
$$
I=-(\tau_p,\tau_s)^oa^{ip}a^{js}
$$
Next, 
$$
II:=\bigl\{(\tau_p,\tau_q)^oa^{ip}\omega_q\bigr\}_{(1)}\tau_j^R=
(\tau_p,\tau_q)^o\langle a^{ip}\omega_q,a^{js}\tau_s\rangle=
(\tau_p,\tau_q)^oa^{ip}a^{jq}
$$
Similarly, 
$$
III:=\tau^R_{i(1)}\bigl\{(\tau_s,\tau_u)^oa^{js}\omega_u\bigr\}=II
$$
and evidently
$$
IV:=\bigl\{(\tau_p,\tau_q)^oa^{ip}\omega_q\bigr\}_{(1)}
\bigl\{(\tau_s,\tau_u)^oa^{js}\omega_u\bigr\}=0
$$
Adding up we get 
$$
j_R(\tau_i)_{(1)}j_R(\tau_j)=(\tau_p,\tau_q)^oa^{ip}a^{jq}
\eqno{(2.9.1)}
$$
Let us differentiate this expression. We have 
$$
\tau_s\{(\tau_p,\tau_q)^oa^{ip}a^{jq}\}=
(\tau_p,\tau_q)^o\bigl\{\tau_s(a^{ip})a^{jq}+a^{ip}\tau_s(a^{jq})\bigr\}=
$$
$$
=(\tau_p,\tau_q)^o\bigl\{-c^{su}_pa^{iu}a^{jq}-c^{sv_q}a^{ip}a^{jv}\bigr\}=
-([\tau_s,\tau_u],\tau_q)^oa^{iu}a^{jq}-
(\tau_p,[\tau_s,\tau_v])^oa^{ip}a^{jv}=0
$$
Therefore, (2.9.1) is a constant. It may be computed by noticing that 
the matrix $(a^{ij})$, considered as a function on the group $G$, 
is equal to the identity at the identity of the group. Hence  
(2.9.1) is equal to $(\tau_i,\tau_j)^o$, which proves (2.5.3). 

{\bf 2.10.} Let us compute 
$$
j_R(\tau_i)_{(0)}j_R(\tau_j)\in \CD^{ch}_{G;(,)1}=T\oplus\Omega
$$
We have
$$ 
j_R(\tau_i)_{(0)}j_R(\tau_j)=j_R(\tau_i)_{(0)}\bigl\{
a^{jq}\tau_q + (\tau_s,\tau_u)^oa^{js}\omega_u\bigr\}
\eqno{(2.10.1)}
$$
Let us compute the first summand. Using (2.8.1), we have  
$$
j_R(\tau_i)_{(0)}(a^{jq}\tau_q)=(j_R(\tau_i)_{(0)}a^{jq})\tau_q=
a^{ip}\tau_p(a^{jq})\tau_q=[\tau_i^R,\tau_j^R]=[\tau_i,\tau_j]^R
\eqno{(9.2)}
$$
by (2.8) and (2.7). 

On the other hand
$$
j_R(\tau_i)_{(0)}\bigl\{(\tau_s,\tau_u)^oa^{js}\omega_u\bigr\}=
\tau^R_{i(0)}\bigl\{(\tau_s,\tau_u)^oa^{js}\omega_u\bigr\}=
(\tau_s,\tau_u)^o\tau_i^R(a^{js}\omega_u)=
$$
by (2.3.2) and (2.2.9)
$$
=(\tau_s,\tau_u)^oa^{ip}\tau_p(a^{js})\omega_u=
(\tau_s,\tau_u)^oc^{ij}_qa^{qs}\omega_u
\eqno{(2.10.3)}
$$
Adding up (2.10.2) and (2.10.3) we see that 
$$
j_R(\tau_i)_{(0)}j_R(\tau_j)=j_R([\tau_i,\tau_j])
$$
which proves (2.5.4) and part (ii) of the theorem. $\btu$ 

{\bf 11.} {\it Proof of} (2.6). The first claim is a reformulation 
of (2.5.3) and (2.5.4).   

The second claim is a trivial consequence 
of (2.5.2) and two Borcherds' 
formulas
$$
x_{(n)}y_{(-1)}z=y_{(-1)}x_{(n)}z+\sum_{j=0}^n\ \binom{n}{j}
(x_{(j)}y)_{(n-1-j)}z
$$
($n\geq 0$), cf. [GII] (0.5.12), and 
$$
(x_{(-1)}y)_{(n)}z=\sum_{j\geq 0}\ \bigl\{
x_{(-1-j)}y_{(n+j)}z+y_{(-n-1-j)}x_{(j)}z\bigr\}
$$
cf. [GII] (0.5.4).  $\btu$

\bigskip\bigskip   
  
\newpage

\centerline{\bf \S 3. BRST} 

\bigskip\bigskip

{\bf 3.1.} Recall the definition of the BRST reduction due to Feigin.  
See [F]; the definition in the language of vertex algebras was given in 
[FF4] (cf. Section 4), cf. also [FF3], Appendix A;  
for a more modern treatment see [BD1], 3.7, [BD2], 7.13. 

Let $\fa$ be a finite dimensional Lie algebra. 
Choose a base $\{a_i\}$ in $\fa$; denote the structure constants 
$$
[a_i,a_j]=c^{ij}_pa_p
\eqno{(3.1.1)}
$$
Recall that the Killing form $( , )_{(K)}:\ \fa\times\fa\lra\BC$ is given 
by 
$$
(a_i,a_j)_{(K)}=c^{ip}_qc^{jq}_p
\eqno{(3.1.2)}
$$
Let $\Pi\fa$ be the space $\fa$ with the reversed parity; denote by 
$\{\phi_i=\Pi a_i\}$ the corresponding base and by $\{\phi^*_j\}$ 
the dual base of $\Pi\fa^*$ given by 
$$
\langle\phi_i,\phi^*_j\rangle=\delta_{ij}
\eqno{(3.1.3)}
$$
Let $C_{BRST}(L\fa)$ denote a graded vertex superalgebra generated by odd 
fields $\phi_i(z)$ of conformal dimension $1$ and odd fields $\phi^*_i(z)$ 
of conformal dimension $0$ with OPE 
$$
\phi_i(z)\phi^*_j(w)\sim\frac{\delta_{ij}}{z-w}
\eqno{(3.1.4)}
$$
We shall identify the spaces $\Pi\fa$ and $\Pi\fa^*$ with their 
obvious images in $C_{BRST}(L\fa)_1$ and $C_{BRST}(L\fa)_0$ 
respectively. 

Let us introduce an odd element $D_{\fa}\in C_{BRST}(L\fa)$ 
of conformal dimension $1$ by 
$$
D_{\fa}=-\frac{1}{2}c^{ij}_p\phi_p\phi^*_i\phi^*_j
\eqno{(3.1.5)}
$$
Thus, we have the corresponding field $D_{\fa}(z)=\sum D_{\fa;n}z^{-n-1}$ 
and we set 
$$
d_{\fa}:=D_{\fa;0}
\eqno{(3.1.6)}
$$
The pair $(C_{BRST}(L\fa),d_{\fa})$ may be regarded as a chiral 
analogue of the Chevalley cochain complex $C(\fa)$. However, in the 
chiral case the square $d^2_{\fa}$ may be nonzero. It is easy to  
compute it. Namely, let us write down the OPE $D_{\fa}(z)D_{\fa}(w)$ 
using Wick theorem. We have 
$$
D_{\fa}(z)D_{\fa}(w)\sim\frac{(a_i,a_j)_{(K)}\phi_i^*(z)\phi^*_j(w)}
{(z-w)^2}
\eqno{(3.1.7)}
$$
Therefore
$$
d_{\fa}^2=(a_i,a_j)_{(K)}\int\phi^*_i(w)'\phi_j(w)
\eqno{(3.1.8)}
$$

{\bf 3.2. Corollary.} {\it If the Lie algebra $\fa$ is nilpotent then 
$d^2_{\fa}=0$.} 

Indeed, the Killing form of a nilpotent Lie algebra is zero. 

{\bf 3.3.} Let $(,):\ \fa\times\fa\lra\BC$ be an arbitrary symmetric invariant 
bilinear form ("level"). Recall that the vertex algebra $\CV_{\fa;(,)}$ 
is generated by 
even fields $a_i(z)$ of conformal weight $1$, subject to OPE 
$$
a_i(z)a_j(w)\sim\frac{(a_i,a_j)}{(z-w)^2}+\frac{[a_i,a_j](w)}{z-w}
\eqno{(3.3.1)}
$$

{\bf 3.4. Lemma.} {\it The rule 
$$
a_i\mapsto c^{ip}_q\phi_q\phi^*_p
\eqno{(3.4.1)}
$$
defines an embedding of vertex algebras 
$$
\CV_{\fa;(,)_{\fa;(K)}}\hra C_{BRST}(L\fa) 
\eqno{(3.4.2)}
$$}  

{\bf 3.5.} Let $\CM$ be a vertex module over $\CV_{\fa; -(,)_{\fa; (K)}}$.   
Let us introduce a space 
$$
C_{BRST}(L\fa;\CM):=C_{BRST}(L\fa)\otimes\CM
\eqno{(3.5.1)}
$$
According to the previous lemma this space is canonically a graded 
(by conformal weight)  
$\CV_{\fa;0}$-supermodule. 
This space is also graded by "fermionic charge" 
$$
C_{BRST}(L\fa;\CM)=\oplus_{p\in\BZ}\ C_{BRST}^{p}(L\fa;\CM)
\eqno{(3.5.2)}
$$
where we assign to $\phi_i$ (resp. $\phi_i^*,\ m\in\CM$) the 
charge $-1$ (resp., $1, 0$). 

Introduce an odd element $D_{\fa;\CM}\in C_{BRST}(L\fa;\CM)$ of 
conformal weight $1$ and fermionic charge $1$ by 
$$
D_{\fa;\CM}=\phi^*_i\otimes a_i+D_{\fa}\otimes 1
\eqno{(3.5.3)}
$$
It follows from (3.1.7) that 
$$
D_{\fa;\CM}(z)D_{\fa;\CM}(w)\sim 0
\eqno{(3.5.4)}
$$
Therefore, setting 
$$
d_{\fa;\CM}:=\int D_{\fa;\CM}(z)
\eqno{(3.5.5)}
$$
we get a differential 
$$
d_{\fa;\CM}^2=0
\eqno{(3.5.6)}
$$
By definition $d$ increases the fermionic charge by $1$. 
 
The pair $(C_{BRST}(L\fa;\CM),d_{\fa;\CM})$ is called the 
{\it BRST complex} of $L\fa$ 
with coefficients in $\CM$, and its cohomology 
$H_{BRST}^{*}(L\fa;\CM)$ is called the BRST cohomology. 

{\bf 3.6. Example.} Let $N$ be a unipotent algebraic group with the 
Lie algebra $\fn$. 
Consider a $\CV_{\fn;0}$-module $\CD^{ch}_{N;0}$ (note that according to 
1.6 this algebra represents a unique isomorphism class of cdo's over $N$). 
Inside the loop algebra $L\fn=\fn[T,T^{-1}]$ consider two Lie subalgebras: 
$\fn_-$ and $\fn_+$, generated by all elements $\tau T^n\ (\tau\in\fn)$ 
with $n<0$ and 
$n\geq 0$ respectively. Then $\CD^{ch}_{N;0}$ is a free $\fn_-$-module 
and a cofree (i.e. the dual module is free) $\fn_+$-module. 

It follows that 
$$
H_{BRST}^i(L\fn; \CD^{ch}_{N;0})=0\ (i\neq 0);\ H_{BRST}^0(L\fn; 
\CD^{ch}_{N;0})=\BC
\eqno{(3.6.1)}
$$    
cf. [FF3], p. 178.

\bigskip\bigskip

\newpage

\centerline{\bf \S 4. Homogeneous spaces}

\bigskip\bigskip

{\bf 4.1.} Let $X$ be a smooth variety. We have an exact 
triangle 
$$
\Omega_X^{[2,3\rangle}\lra\Omega_X^{[2}\lra C
\eqno{(4.1.1)}
$$
where 
$$
\Omega_X^{[2}:\ 0\lra\Omega^2_X\lra\Omega^3_X\lra\ldots
\eqno{(4.1.2)}
$$
($\Omega^2_X$ sitting in degree $0$) and 
$$
C:\ 0\lra\Omega^4_X/d\Omega^3_X\lra\Omega^5_X\lra\ldots
\eqno{(4.1.3)}
$$
($\Omega^4_X/d\Omega^3_X$ sitting in degree $2$). It follows that 

{\bf 4.1.1.}\footnote{this remark is due to H.Esnault} 
{\it the canonical map 
$$
H^i(X;\Omega_X^{[2,3\rangle})\lra H^i(X,\Omega_X^{[2})
\eqno{(4.1.4)}
$$ 
is injective for $i=2$ and bijective for $i=0,1$.}  

{\bf 4.2.} Let $G$ be a simple  algebraic group. 
In this section we shall discuss the chiral differential operators 
on homogeneous spaces $G/G'$ where $G'=N$ --- a unipotent subgroup, 
$G'=P$ --- a parabolic but not minimal parabolic, or 
$G'=B$ --- a Borel subgroup. 

\bigskip

{\it The case $G/N$}

\bigskip

{\bf 4.3.} Consider the projection $\pi:\ G\lra X:=G/N$. The variety 
$X$ is quasiaffine, therefore we have 
$$
H^i(X;\Omega_X^{[2})=H^i\Gamma(X;\Omega_X^{[2})=
H^{i+2}_{DR}(X)\ (i\geq 1)
\eqno{(4.3.1)}
$$
On the other hand, 
$$
H^*_{DR}(X)=H^*(X;\BC)
\eqno{(4.3.2)}
$$
by Grothendieck's theorem, (cf. [Gr], Theorem 1'). 

The projection $\pi$ is an affine morphism which is a Zariski locally 
trivial bundle with fiber $N$ 
isomorphic to an affine space, so $\pi^*:\ H^*(X;\BC)\iso H^*(G;\BC)$. 
It follows that 
$$
\pi^*:\ H^i(X;\Omega_X^{[2})\iso H^i(G;\Omega_G^{[2})
\ (i\geq 1)
\eqno{(4.3.3)}
$$   
We have a short exact sequence 
$$
0\lra \Theta_{G/X}\lra \Theta_G\lra\pi^*\Theta_X\lra 0
\eqno{(4.3.4)}
$$
and the vector bundles $\Theta_G, \Theta_{G/X}$ are trivial 
(a base of global sections of $\Theta_{G/X}$ is given by left 
invariant vector fields coming from the Lie algebra $\fn:=Lie(N)$). 

Therefore we have 
$$
\pi^*c(\fD iff^{ch}_X)=\pi^*ch_2(\Theta_X)=ch_2(\Theta_G)=0
\eqno{(4.3.5)}
$$
hence
$$
c(\fD iff^{ch}_X)=0
\eqno{(4.3.6)}
$$
by 4.1.1. 

Therefore, 

{\bf 4.4.} {\it the groupoid $\Gamma(X;\fD iff^{ch}_X)$ is nonempty.  
The set of isomorphism classes $\pi_0(\Gamma(X,\fD iff^{ch}_X))$ is 
a torseur under $H^3_{DR}(X)=H^3_{DR}(G)$.} 

{\bf 4.5.} Let $\CD^{ch}_{G;(,)}$ be the sheaf of chiral differential 
operators on $G$ of level $(,)$ where 
$$
(,):\ \fg\times\fg\lra\BC
\eqno{(4.5.1)}
$$
is a fixed symmetric invariant bilinear form on $\fg=Lie(G)$, cf. 1.7.  
The form (4.5.1) is a scalar multiple of the Killing form 
$$
(,)=c(,)_{\fg;(K)},\ c\in\BC
\eqno{(4.5.2)}
$$
The Killing form on $\fg$ restricts to zero on $\fn$ (since the trace 
of a nilpotent 
endomorphism is zero). Therefore we have the canonical embedding 
of vertex algebras 
$$
\CV_{\fn;0}\hra\CV_{G;(,)}\hra\CD^{ch}_{G;(,)}
\eqno{(4.5.3)}
$$
so that the sheaf $\CD^{ch}_{G;(,)}$ becomes a sheaf of $\CV_{\fn;0}$-modules. 

Applying the BRST construction 3.5 to $\CM=\pi_*\CD^{ch}_{G;(,)}$ and 
$\fa=\fn$ we get a sheaf of BRST complexes 
$C_{BRST}(L\fn;\pi_*\CD^{ch}_{G;(,)})$ and BRST cohomology sheaves  
$H_{BRST}^{*}(L\fn;\pi_*\CD^{ch}_{G;(,)})$ over $X$. They are sheaves 
of $\BZ_{\geq 0}$-graded vertex superalgebras. 

{\bf 4.6. Theorem.} {\it We have 
$$
H_{BRST}^{i}(L\fn;\pi_*\CD^{ch}_{G;(,)})=0\ (i\neq 0)
\eqno{(4.6.1)}
$$
The sheaf $H_{BRST}^{0}(L\fn;\pi_*\CD^{ch}_{G;(,)})$ is an algebra of chiral 
differential operators over $X$. 

The correspondence 
$$
\CD^{ch}_{G;(,)}\mapsto H_{BRST}^{0}(L\fn;\pi_*\CD^{ch}_{G;(,)})
\eqno{(4.6.2)}
$$
induces a bijection of the sets of isomorphism classes 
$$
\pi_0(\Gamma(G;\fD iff^{ch}_G))\iso 
\pi_0(\Gamma(X;\fD iff^{ch}_X))
\eqno{(4.6.3)}
$$}

We shall use the notation $\CD^{ch}_{X;(,)}$ for the cdo 
$H_{BRST}^{0}(L\fn;\pi_*\CD^{ch}_{G;(,)})$.  

Note that higher direct images 
$R^i\pi_*\CD^{ch}_{G;(,)}$ are trivial for $i>0$ since the morphism $\pi$ 
is affine and the sheaves $\CD^{ch}_{G;(,)}$ admits a filtration 
whose qutients are coherent (in fact, locally free) $\CO_G$-modules. 

{\bf 4.7.} {\it Proof (sketch).} Locally on $X$ the projection 
$\pi: G\lra X$ is isomorphic to the direct product $U\times N\lra U$. 
If $\CD$ is an algebra of chiral differential operators on $U\times N$ then 
$\CD\iso \CD_U\boxtimes\CD_N$ where $\CD_U$ (resp. $\CD_N$ is an 
algebra of differential operators on $U$ (resp. $N$). 
Now the first claim of the theorem follows from (3.6.1).  

The second claim is a corollary of 4.4. $\btu$ 

{\bf 4.8. Corollary.} {\it The dual embedding   
$$
j_R:\ \CV_{G; (,)^o}\lra \CD^{ch}_{G; (,)}
\eqno{(4.8.1)}
$$
defined in} Corollary 2.6 {\it induces a canonical morphism of vertex 
algebras 
$$
i^o_{(,)}:\ \CV_{G; (,)^o}\lra \CD^{ch}_{X;(,)}
\eqno{(4.8.2)}
$$} 

Indeed, we know that (4.8.1) commutes with the left action of $\hfg$ 
hence with the BRST differential. 

In particular, the Kac-Moody algebra $\hfg$ at level $(,)^o$ 
acts canonically on {\it Wakimoto modules} which may be defined as  
the spaces of sections $\Gamma(U;\CD^{ch}_{X,(,)})$ where $U$ is 
a big cell. This is a result due to Feigin-Frenkel obtained in [FF1] - [FF3]   
in a different way.

\bigskip

{\it The case $G/B$}

\bigskip

{\bf 4.9.} Let $B\subset G$ be a Borel subgroup, $\pi:\ G\lra X:=G/B$. 
$X$ is a smooth projective variety and  
we have 
$$
H^p(X;\Omega^q_X)=0\ (p\neq q)
\eqno{(4.9.1)}
$$
and 
$$
H^i(X;\Omega^i_X)=H^{2i}(X;\BC)
\eqno{(4.9.2)}
$$
It follows that $H^2(X;\Omega^{[2})=H^4(X;\BC)$, and [GII], Theorem 7.5 
says that the image of $c(\fD iff^{ch}_X)$ in $H^4(X;\BC)$ is equal 
to
$$
2ch_2(\Theta_X):=c_1^2(\Theta_X)-c_2(\Theta_X)
\eqno{(4.9.3)}
$$
where $c_i(\Theta_X)\in H^{2i}(X;\BC)$ are the Chern classes of the tangent 
bundle $\Theta_X$. 

{\bf 4.9.1. Lemma.} {\it We have 
$$
ch_2(\Theta_X)=0
\eqno{(4.9.1.1)}
$$}

{\it Proof.} The space $H^2(X;\BC)$ may be identified with the complexification 
of the root lattice of $G$. The classical theorem by J. Leray says that 
that the cohomology algebra $H^*(X;\BC)$ is equal to the quotient of 
the symmetric algebra of the space $H^2(X;\BC)$ modulo the ideal 
generated by the subspace of invariants of the Weyl group $W$ having positive 
degree, cf. [L], Th\'eoreme 2.1 b).  

The class $[\Theta_X]$ in the Grothendieck 
group of vector bundles is equal to the sum $\sum\ [\CL_{\alpha}]$ where 
$\alpha$ runs through all negative roots, and $\CL_{\alpha}$ is 
the line bundle with $c_1(\CL_{\alpha})=\alpha$. Hence $ch_2(\Theta_X)=
\sum\ \alpha^2$; this element is invariant under the action 
of $W$, therefore its image in $H^4(X;\BC)$ is zero. $\btu$       
          
This lemma implies that  
$$
c(\fD iff^{ch}_X)=0
\eqno{(4.9.4)}
$$
by 4.1.1. On the other hand, it follows from 
4.1.1 and (4.9.1) that 
$$
H^i(X;\Omega_X^{[2,3\rangle})=0\ (i=0,1)
\eqno{(4.9.5)}
$$
Thus, by [GII], {\it loc. cit.} we get 

{\bf 4.10. Theorem.} {\it The groupoid $\Gamma(X;\fD iff^{ch}_X)$ is nonempty 
and trivial. In other words, there exists a unique, up to a unique 
isomorphism, algebra of chiral differential operators $\CD^{ch}_X$ over $X$.} 
$\btu$

{\bf 4.11.} Let us construct the algebra $\CD_X^{ch}$ using BRST reduction. 
Consider the sheaf $\CD^{ch}_{G;(,)}$ as in 4.5. Let 
$(,)_{\fb}$ denote the restriction of $(,)$ to $\fb:=Lie(B)$. 

We have a canonical embedding of vertex algebras 
$$
\CV_{\fb;(,)_{\fb}}\hra \CD^{ch}_{G;(,)}
\eqno{(4.11.1)}
$$
We have 
$$
(,)_{\fg; (K)}|_{\fb\times\fb}=2(,)_{\fb;(K)}
\eqno{(4.11.2)}
$$
Therefore, to get the minus Killing on $\fb$ we 
have to start from $-(,)_{\fg, (K)}/2$, i.e. from the 
critical level on $\fg$. Thus, by construction 3.5,  

{\bf 4.11.1.} {\it the BRST complex $C_{BRST}(L\fb;\pi_*\CD^{ch}_{G;(,)})$ 
is defined iff $(,)=(,)_{\fg;crit}$, i.e. on the} {\bf critical 
level}, cf. (2.1.5).  

Let us denote the algebra $\CD^{ch}_{G;(,)_{\fg; crit}}$ by 
$\CD^{ch}_{G;crit}$. 

{\bf 4.12.} Let $\CV_{\fb; crit}$ denote the vacuum module on the 
critical level $\CV_{\fb;(,)_{\fb; crit}}$. Let $\CM$ be a vertex 
algebra equipped with a morphism of vertex algebras  
$\CV_{\fb; crit}\lra\CM$. 

We have a decomposition 
$$
\fb=\fh\oplus\fn
\eqno{(4.12.1)}
$$
where $\fn\subset\fb$ (resp. $\fh\subset\fb$) the maximal nilpotent 
(resp. the Cartan) subalgebra; $\fn$ is spanned by  elements $e_{\alpha}$, 
$\alpha$ 
being a positive root. Let us define a 
vertex subalgebra 
$$
C_{BRST}(L\fb,\fh;\CM)\subset C_{BRST}(L\fb;\CM)
\eqno{(4.12.2)}
$$ 
The vector space $C_{BRST}(L\fb;\CM)$ is spanned by  
all the monomials 
$$
h_{1;i_1}\cdot\ldots \cdot h^*_{1;j_1}\cdot\ldots\cdot 
e_{\alpha_1;k_1}\cdot\ldots\cdot 
e^*_{\alpha'_1;l_1}\cdot\ldots\otimes m_{\mu}
\eqno{(4.12.3)}
$$     
where 
the indices $i_p, j_p, k_p, l_p$ denote the conformal weight, \
$i_p, k_p\geq 1;\ j_p, l_p\geq 0$,\ $m_{\mu}\in\CM$ is a vector 
of weight $\mu\in\fh^*$. 

By definition, the subspace $C_{BRST}(L\fb,\fh;\CM)$ is spanned by all 
monomials (4.12.3) such that 

(a) all $j_p\geq 1$; 

(b) $\sum \alpha_p -\sum \alpha'_q+\mu=0$. 

It is a vertex subalgebra. One checks that the BRST differential $d$ 
in $C_{BRST}(L\fb;\CM)$ preserves $C_{BRST}(L\fb,\fh;\CM)$. 

We define the relative BRST cohomology $H_{BRST}^{*}(L\fb,\fh;\CM)$ as 
the cohomology of $C_{BRST}(L\fb,\fh;\CM)$ with respect to this 
differential. It is canonically a vertex algebra. 

{\bf 4.13.} Applying the previous definition to $\CM=\pi_*\CD^{ch}_{G;crit}$ 
we get the BRST cohomology sheaves 
$H_{BRST}^{*}(L\fb,\fh;\pi_*\CD^{ch}_{G;crit})$ on $X$. 

{\bf 4.14. Theorem.} {\it We have 
$$
H_{BRST}^{i}(L\fb,\fh;\pi_*\CD^{ch}_{G;crit})=0\ (i\neq 0)
\eqno{(4.14.1)}
$$
and the sheaf of vertex algebras 
$H_{BRST}^{0}(L\fb,\fh;\pi_*\CD^{ch}_{G;crit})$ 
is canonically isomorphic to $\CD^{ch}_X$.} 

The proof is the same as that of Theorem 4.6. 

{\bf 4.15. Corollary.} {\it We have a canonical isomorphism 
of vertex algebras
$$
H^*(X;\CD^{ch}_X)=H_{BRST}^{*}(L\fb,\fh;\Gamma(G;\CD^{ch}_{G;crit}))
\eqno{(4.15.1)}
$$}    

Indeed, as we have already remarked, sheaves of cdo $\CD^{ch}_G$ on $G$ 
have a filtration whose graded quotients are vector bundles; 
but the vatiety $G$ and the morphism $\pi$ are affine, hence 
$H^i(G;\CD^{ch}_G)=R^i\pi_*\CD^{ch}_G=0$ for $i>0$. 

{\bf 4.16. Corollary.} {\it The dual embedding   
$$
j_R:\ \CV_{G; crit}\lra \CD^{ch}_{G; crit}
\eqno{(4.16.1)}
$$
defined in} Corollary 2.6 {\it induces a canonical morphism of vertex 
algebras 
$$
i_{crit}:\ \CV_{G; crit}\lra \CD^{ch}_X
\eqno{(4.16.2)}
$$} 

Indeed, we know that (4.16.1) commutes with the left action of $\hfg$ 
hence with the BRST differential. 

In particular, the Kac-Moody algebra $\hfg$ on the critical level 
acts canonically on {\it Wakimoto modules} which may be defined as  
the spaces of sections $\Gamma(U;\CD^{ch}_X)$ where $U$ is 
a big cell. This is a result of Feigin-Frenkel proven in [FF1] - [FF3]   
in a different way. 

\bigskip

{\it The case $G/P$} 

\bigskip

{\bf 4.17.} Let $P\subset G$ be a parabolic but not Borel, $\fp=Lie(P)$,  
$\pi:\ G\lra X:=G/P$. The discussion of 4.9 applies as it is, except for 
Lemma 4.9.1, which is replaced by 

{\bf 4.17.1. Lemma.} {\it We have 
$$
ch_2(\Theta_X)\neq 0
\eqno{(4.17.1.1)}
$$}  

Hence $c(\fD iff^{ch}_X)\neq 0$, and we get 

{\bf 4.18. Theorem.} {\it The groupoid $\Gamma(X;\fD iff^{ch}_X)$ is empty. 
In other words, there is no cdo over $X$.} $\btu$ 

{\bf 4.19.} Let us see how this fact is reflected in the BRST world. 
We would like to get a sheaf of cdo over $X$ as the BRST cohomology 
of $L\fp$ with coefficients in some sheaf $\pi_*\CD^{ch}_{G;(,)}$. 
However, no form $(,)$ on $G$ restricts to the Killing form on $\fp$, 
which implies that the square of the BRST differential in 
$C_{BRST}(L\fp;\pi_*\CD^{ch}_{G;(,)})$ is {\it always} nonzero. 
Thus, the BRST cohomology is not defined, which is compatible with 4.18.

\bigskip\bigskip

%\newpage 

%\input groupref
\centerline{\bf References}

\bigskip\bigskip

[BD1] A.~Beilinson, V.~Drinfeld, Chiral algebras, Preprint. 

[BD2] A.~Beilinson, V.Drinfeld, Quantization of Hitchin integrable 
system and Hecke eigensheaves, Preprint. 

%[B] A.~Borel, Topology of Lie groups and characteristic classes, 
%{\it Bull. Amer. Math. Soc.}, {\bf 61} (1955), 397-432. 

%[BH] A.~Borel, F.~Hirzebruch, Characteristic classes and homogeneous 
%spaces, I, {\it Amer. J. Math.}, {\bf 80}, 458-538. 

[F] B.L.~Feigin, Semiinfinite homology of Lie, Kac-Moody and 
Virasoro algebras, {\it Uspekhi Matem. Nauk}, {\bf 39}, no. 2(236) 
(1984), 195-196 (Russian). 

[FF1] B.L.~Feigin, E.~Frenkel, A family of representations of affine Lie 
algebras, {\it Uspekhi Matem. Nauk}, {\bf 43}, no. 5(263) (1988), 227-228 
(Russian). 

[FF2] B.L.~Feigin, E.~Frenkel, Representations of affine Kac-Moody 
algebras and bosonization, {\it Physics and Mathematics of Strings}, 271-316. 
World Sci. Publishing, Teaneck, NJ, 1990. 

[FF3] B.L.~Feigin, E.~Frenkel, Affine Kac-Moody algebras and semi-infinite 
flag manifolds, {\it Comm. Math. Phys.}, {\bf 128} (1990), 161-189. 

[FF4] B.~Feigin, E.~Frenkel, Affine Kac-Moody algebras at the critical 
level and Gelfand-Dikii algebras, {\it Infinite analysis, Part A, B 
(Kyoto, 1991)}, 197 - 215, {\it Adv. Ser. Math. Phys.}, {\bf 16}, 
World Sci. Publishing, River Edge, NJ, 1992.   
   
[GII] V.~Gorbounov, F.~Malikov, V.~Schechtman, Gerbes of chiral 
differential operators. II, math.AG/0003170. 

[Gr] A.~Grothendieck, On the de Rham cohomology of algebraic 
varieties, {\it Publ. Math. IHES}, {\bf 29} (1966), 95-103. 

[L] J.~Leray, Sur l'homologie des groupes de Lie, des espaces 
homog\`enes et des espaces fibr\'es principaux, {\it Colloque de Topologie 
(Espaces Fibr\'es) (Bruxelles, 1950)}, 101-115, Georges Thone, Li\`ege, 
Masson \& Cie, Paris, 1951.   

[MSV] F.~Malikov, V.~Schechtman, A.~Vaintrob, Chiral de Rham complex, 
{\it Commun. Math. Phys.}, {\bf 204} (1999), 439-473.

\bigskip

\bigskip

V.G.: Department of Mathematics, University of Kentucky, 
Lexington, KY 40506, USA;\ vgorb\@ms.uky.edu

F.M.: Department of Mathematics, University of Southern California, 
Los Angeles, CA 90089, USA;\ fmalikov\@mathj.usc.edu   

V.S.: IHES, 35 Route de Chartres, 91440 Bures-sur-Yvette, France;\ 
vadik\@ihes.fr

\enddocument